\documentclass{article}
\usepackage{amsmath,amsthm,amscd,amssymb,amsfonts}
\input{amssym.def}
\input{amssym.tex}
\newtheorem{df}{Definition}[section]
\newtheorem{thm}{Theorem}[section]
\newtheorem{prop}{Proposition}[section]
\newtheorem{lm}{Lemma}[section]

\newtheorem{remark}{Remark}[section]
\newtheorem{fact}{Fact}[section]
\newtheorem{cor}{Corollary}[section]

\title{A quantization of the Hitchin hamiltonian system and the
Beilinson-Drinfeld isomorphism}
\author{Ken-ichi SUGIYAMA
\footnote{Address : Ken-ichi SUGIYAMA,
 Department of Mathematics and Informatics,
Faculty of Science,
Chiba University,
1-33 Yayoi-cho Inage-ku,
Chiba 263-8522, Japan}
\footnote{e-mail address : sugiyama@math.s.chiba-u.ac.jp}}
\begin{document}
\maketitle
\begin{abstract}Let $X$ be a smooth projective curve defined over ${\mathbb C}$ whose genus is
 greater than one.
We will generalize the Hitchin's hamiltonian system to the cotangent
 bundle of the modular stack of principal $SL_{2}({\mathbb C})$ bundle with
 parabolic reductions on $X$. Also its properties will be investigated. As an
 application, we will show a generalization of the Beilinson-Drinfeld
 isomorphism, which is a quantization of the Hitchin system.  
\end{abstract}
\section{Introduction}
Let $X$ be a connected smooth projective
curve defined over ${\mathbb C}$. Suppose we are given a local system of
rank one ${\cal L}$ on $X$.  Then it determines a representation
\[
 \pi_1(X,\,x_0)\stackrel{\rho_{\cal L}}\longrightarrow {\mathbb C}^{\times},
\]
which factors through the maximal abelian quotient of the fundamental
group $H_1(X,\,{\mathbb Z})$. If we identify $H_1(X,\,{\mathbb
Z})$ with the fundamental group of the Jacobian ${\rm Jac}(X)$,
$\rho_{\cal L}$ defines a local system ${\cal M}_{{\cal L}}$ of rank one
on ${\rm Jac}(X)$. Such a correspondence should be considered as {\it a
geometric abelian class field theory.}
\par\vspace{5pt}
In order to pass to a nonabelian situation, we have to change a point of
view to understand this phenomenon (\cite{Frenkel2003}). The local system ${\cal L}$ defines a flat connection
\[
 \nabla=d\,+\,A,\quad A\in H^0(X,\,\Omega),
\]
where $\Omega$ is the canonical line bundle of $X$. Note that the
cotangent bundle $T^*{\rm Jac}(X)$ is the trivial
vector bundle, which may be identified with
\[
 {\rm Jac}(X)\times H^0(X,\,\Omega). 
\]
Let us observe that the projection
\[
 T^*{\rm Jac}(X)\stackrel{p}\longrightarrow H^0(X,\,\Omega)
\]
induces an isomorphism of the global coordinate rings:
\[
 \Gamma(H^0(X,\,\Omega),\,{\cal O})\stackrel{p^*}\simeq \Gamma(T^*{\rm
 Jac}(X),\,{\cal O}).
\]
On the other hand, one easily see that the ring of global differential operators $D({\rm Jac}(X))$ on the Jacobian is isomorphic to $\Gamma(T^*{\rm
 Jac}(X),\,{\cal O})$ by the symbol map. Combining this with the
 isomorphism above, we obtain
\[
 \Gamma(H^0(X,\,\Omega),\,{\cal O}) \simeq D({\rm Jac}(X)). 
\]
Now the connection form $A$ determines a homomorphism
\[
 D({\rm Jac}(X))\simeq \Gamma(H^0(X,\,\Omega),\,{\cal O})\stackrel{f_A}\longrightarrow{\mathbb C},
\]
which yields a D-module ${\cal M}_A$ on ${\rm Jac}(X)$:
\[
 {\cal M}_A={\cal D}_{{\rm Jac}(X)}\otimes_{D({\rm Jac}(X))}{\mathbb C}.
\]
Here ${\cal D}_{{\rm Jac}(X)}$ is the sheaf of differential operators on
the Jacobian. Then ${\cal M}_A$ corresponds to the local system ${\cal
M}_{{\cal L}}$.
\par\vspace{5pt}
The latter picture has the following nonabelization due to Beilinson and
Drinfeld (\cite{Beilinson-Drinfeld}). Let $G$ be a semisimple Lie group and let ${\rm Bun}_{G,X}$ be
the classifying stack of principal $G$-bundles on $X$. Hitchin
has defined a map which is referred as {\it the Hitchin's hamiltonian
system} (\cite{Hitchin1987}):
\[
 T^*{\rm Bun}_{G,X}\stackrel{h}\longrightarrow H_{G}.
\]
Here $H_G$ is an affine space which depends on $G$. Note that the map
$h$ corresponds to the projection $p$ in the geometric abelian class
field theory. Beilinson and Drinfeld have {\it quantized} the Hitchin's
system. Namely they have constructed an isomorphism between the ring of
global differential operators on a half canonical of ${\rm Bun}_{G,X}$
and the global coordinate ring of $H_G$. Using it, they
have established {\it unramified geometric nonabelian class field
theory}, which assigns a holonomic D-module on ${\rm Bun}_{G,X}$ to an
element of $H_G$.
\par\vspace{5pt}
We will treat a case of which a local system is allowed to be ramified at several
points. 
We will restrict ourselves to the case that $G$ is $SL_2({\mathbb C})$.
\par\vspace{5pt}
Let $\{z_1,\cdots,z_N\}$ be distinct points on $X$ and let us form a
divisor
\[
 D=\sum_{i=1}^N z_i.
\]
We assume that the degree of the line bundle $\Omega(D)$ is positive.
We will consider the classifying stack ${\rm Bun}_{G,X}^{D-fl}$ of
principal $G$-bundles with $D$-flags on $X$ (see $\S 2$). In $\S 3$, for
$\lambda=(\lambda_1,\cdots, \lambda_N)\in {\mathbb Z}^N$, we will construct
a sheaf of twisted differential operators ${\cal
D}^{\prime}_{D-fl,\lambda}$ on ${\rm Bun}_{G,X}^{D-fl}$. Let
$D^{\prime}_{D-fl,\lambda}$ be the ring of 
its global sections. Then we will show the following theorem (see $\S 5$).
\begin{thm}
There is a isomorphism of ${\mathbb C}$-algebras:
\[
 \Gamma(H^0(X,\,\Omega^{\otimes 2}(D)),\,{\cal O})\stackrel{h_{\lambda}}\simeq D^{\prime}_{D-fl,\lambda}. 
\]
\end{thm}
When $X$ is the projective line, {\bf Theorem 1.1} has been already
established by Frenkel (\cite{Frenkel1995}). As before, this will
establish a geometric {\it ramified} nonabelian class field theory in
the case of $SL_2({\mathbb C})$, which assigns a holonomic D-module on ${\rm
Bun}_{G,X}^{D-fl}$ to an element of $H^0(X,\,\Omega^{\otimes 2}(D))$
(see $\S 5$). It
is conjectured that the D-module should be regular holonomic and moreover
that it should be a Hecke eigensheaf. We will treat the conjecture in
the near furure.
\par\vspace{5pt}
Throughout the paper, we set $G=SL_2({\mathbb C})$ and its Lie algebra will
be denoted by  $\frak{G}$.
\section{The Hitchin system for principal $G$ bundles with $D$-flags}
Let $X$ be a connected smooth projective curve defined over ${\mathbb C}$
of genus $g$. 
Let $\{z_1,\cdots,z_N\}$ be mutually distinct points of $X$ and we set
\[
 D=\sum_{i=1}^{N}z_i \in {\rm Div}(X).
\]
We will always assume that the degree of both of $\Omega(D)$ and ${\cal O}(D)$ are positive. Here $\Omega$ is the canonival line bundle of $X$.
We will choose and fix a local coordinate $t_i$ at each $z_i$. In this
section, we will
study basic properties of the modular stack of principal $G$ bundles with $D$-flags. See \cite{LS}, \cite{Laumon-MoretBailly} or \cite{Sorger} for a treatment of stacks.
\subsection{The modular stack of principal $G$ bundles with $D$-flags}
Let $P$ be a principal $G$-bundle on $X$. {\it A $D$-flag} of $P$ is
defined to be an $N$-tuple 
\[
 \{B_1,\cdots,B_N\},
\]
where $B_i$ is a Borel subgroup of $P|_{z_i}$. Let $Aff_{\slash{\mathbb
C}, et}$ be the category of affine schemes defined over ${\mathbb C}$ with
the \'{e}tale topology. We will consider the modular stack
\[
 {\rm Bun}_{G,X}^{D-fl}
\]
over $Aff_{\slash{\mathbb C}, et}$ which classifies principal $G$-bundles on
$X$ endowed with $D$-flags. It is a smooth stack of dimension
$3(g-1)+N$ and has the following structure. Note that since we have
assumed that degree of $\Omega(D)$ is positive, $3(g-1)+N$ is nonnegative.
\par\vspace{5pt}
Let ${\rm Bun}_{G,X}$ be the modular stack on $Aff_{\slash{\mathbb C}, et}$
which classifies principal $G$-bundles on $X$. Forgetting $D$-flags, we
have a morphism
\[
 {\rm Bun}_{G,X}^{D-fl}\stackrel{\pi}\longrightarrow {\rm Bun}_{G,X},
\]
which makes ${\rm Bun}_{G,X}^{D-fl}$ into a $({\mathbb P}^1)^N$-fibration
over ${\rm Bun}_{G,X}$. Since ${\rm Bun}_{G,X}$ is connected and smooth, ${\rm Bun}_{G,X}^{D-fl}$ is also connected. 
\subsection{The cotangent bundle of ${\rm Bun}^{D-fl}_{G,x}$}
We will study the cotangent space of ${\rm Bun}_{G,X}^{D-fl}$ at
$(P,\{B_i\}_i)$.
Let ${\frak B}_i$ be the Lie algebra of $B_i$ and we set
${\frak N}_i={\frak G}\slash{\frak B}_i$, which will be regarded as a nilpotent subalgebra of ${\frak G}$.
We will consider a skyscraper sheaf supported on $Supp(D)$:
\[
 {\frak G}_{D}={\frak G}^{\oplus N},
\] 
and its subsheaves
\[
 {\frak B}_{D}=\oplus_{i=1}^{N}{\frak B}_i,\quad {\frak N}_{D}=\oplus_{i=1}^{N}{\frak N}_i.
\]
Let ${\rm ad}_P({\frak G})$ be the adjoint ${\frak G}$-bundle associated
to $P$. A sheaf ${\frak B}_{P,D}$ is defined to be
\[
 {\frak B}_{P,D}={\rm Ker}[{\rm ad}_P({\frak
 G})\stackrel{ev_D}\longrightarrow {\frak N}_{D}].
\]
Here $ev_D$ is the composition of the evaluation map at $D$:
\[
 {\rm ad}_P({\frak G}) \longrightarrow {\frak G}_{D}
\]
with the natural projection onto ${\frak N}_{D}$. Note that we have the
following commutative diagram:
\begin{equation}
\begin{array}{ccccccccc}
0 & \to & {\rm ad}_{P}({\frak G})(-D) & \to & {\rm ad}_{P}({\frak G}) & \to & {\frak G}_{D} & \to &0\\
  &     & \downarrow &  & \parallel & & \downarrow & & \\
0 & \to & {\frak B}_{P,D} & \to & {\rm ad}_{P}({\frak G}) & \to & {\frak N}_{D} & \to &0.
\end{array}
\end{equation}
such that the horizontal sequences are exact.
By the deformation theory, the tangent space of ${\rm Bun}_{G,X}^{D-fl}$
at $(P,\{B_i\}_i)$ is isomorphic to $H^{1}(X,\,{\frak
B}_{P,D})$. Therefore by the Serre
duality, the cotangent space is identified with
\[
 H^{0}(X,\,{Hom}_{{\cal O}_X}({\frak B}_{P,D},{\cal O}_X)\otimes_{{\cal O}_X} \Omega_X).
\]
We want to rewrite this more convenient form.
\par\vspace{5pt}
By the snake lamma, (1) shows
\[
 {\rm Coker}[{\rm ad}_{P}({\frak G})(-D) \longrightarrow {\frak B}_{P,D}]
\]
is isomorphic to ${\frak B}_D$, which implies an exact seqence
\begin{equation}
 0 \to {\rm ad}_{P}({\frak G})(-D) \to {\frak B}_{P,D} \to {\frak B}_D
 \to 0.
\end{equation}

By the Killing form, we will identify ${Hom}_{{\cal O}_X}({\rm
ad}_P({\frak G})(-D),{\cal O}_X)$ with ${\rm ad}_{P}({\frak G})(D)$. Then
(2) implies
\[
 0 \to {Hom}_{{\cal O}_X}({\frak B}_{P,D},{\cal O}_X)\otimes_{{\cal
 O}_X} \Omega_X \to {\rm ad}_{P}({\frak G})(D)\otimes_{{\cal
 O}_X}\Omega_X \to 
Ext^{1}_{{\cal O}_X}({\frak B}_D,{\cal O}_X)\otimes_{{\cal
 O}_X} \Omega_X \to 0.
\] 
The residue map shows us that the last term is isomorphic to ${\frak B}_D$ and
therefore we have
\[
 0 \to {Hom}_{{\cal O}_X}({\frak B}_{P,D},{\cal O}_X)\otimes_{{\cal
 O}_X} \Omega_X \to {\rm ad}_{P}({\frak G})(D)\otimes_{{\cal
 O}_X}\Omega_X \to {\frak B}_D \to 0.
\]
Thus we have proved the following proposition.
\begin{prop}
The cotangent space $T^{*}_{(P,\{B_i\}_i)}{\rm Bun}_{G,X}^{D-fl}$ 
is isomorphic to
\[
 {\rm Ker}[H^{0}(X,\,{\rm ad}_{P}({\frak G})(D)\otimes_{{\cal
 O}_X}\Omega_X)\stackrel{\rm Res}\longrightarrow {\frak B}_D].
\]

\end{prop}

\begin{remark} Using the snake lemma, the commutative diagram whose the middle and all horizontal lines are exact:
\[
\begin{array}{ccccccccc}
 &  &  &  & 0 &  &  &  &  \\
 &  &  &  & \downarrow &  &  &  &  \\
 &  &  &  & {\rm ad}_{P}({\frak G})\otimes_{{\cal O}_X} \Omega_{X} &  &  &  & \\
  &     &  &  & \downarrow & &  & & \\
0 & \to & Hom_{{\cal O}_X}({\frak B}_{P,D},{\cal O}_X)\otimes_{{\cal O}_X} \Omega_X & \to & {\rm ad}_{P}({\frak G})(D)\otimes_{{\cal O}_X} \Omega_{X} & \stackrel{{\rm Res}}\to & {\frak B}_{D}\otimes_{{\cal O}_X}\Omega_X & \to &0\\
  &     & \downarrow &  & \downarrow & & \parallel & & \\
0 & \to & {\frak N}_{D}\otimes_{{\cal O}_X}\Omega_X & \to & {\frak G}_{D}\otimes_{{\cal O}_X}\Omega_X & \to & {\frak B}_{D}\otimes_{{\cal O}_X}\Omega_X & \to &0 \\
 &  &  &  & \downarrow &  &  &  &  \\
&  &  &  & 0 &  &  &  &  
\end{array}  
\]
shows that there is an exact sequence
\[
 0 \to {\rm ad}_{P}({\frak G})\otimes_{{\cal O}_X} \Omega_{X} \to Hom_{{\cal O}_X}({\frak B}_{P,D},{\cal O}_X)\otimes_{{\cal O}_X} \Omega_X \to {\frak N}_{D}\otimes_{{\cal O}_X}\Omega_X  \to 0.
\]
\end{remark}

\subsection{The Hitchin's hamiltonian system}
{\bf Proposition 2.1} shows $T^{*}_{(P,\{B_i\}_i)}{\rm
Bun}_{G,X}^{D-fl}$ may be identified a subspace of 
\[
 H^{0}(X,\,{\rm ad}_{P}({\frak G})(D)\otimes_{{\cal
 O}_X}\Omega_X)
\]
which consists of $A$ which has a Taylor expansion at each $z_i$:
\[
 A=
\left(
\begin{array}{cc}
a_i(t_i) & b_i(t_i)\\
c_i(t_i) & -a_i(t_i)
\end{array}
\right)
\frac{dt_i}{t_i}, \quad a_i(t_i), b_i(t_i), c_i(t_i)\in {\mathbb C}[[t_i]], 
\]
satisfying 
\[
\left(
\begin{array}{cc}
a_i(0) & b_i(0)\\
c_i(0) & -a_i(0)
\end{array}
\right)
\in {\frak N}_i. 
\] 
In particular for $A\in T^{*}_{(P,\{B_i\}_i)}{\rm
Bun}_{G,X}^{D-fl}$, we have
\[
 \det A \in H^0(X,\,\Omega^{\otimes 2}(D)),
\]
and when $(P,\{B_i\}_i)$ moves over ${\rm
Bun}_{G,X}^{D-fl}$, the determinant defines a morphism which will be
referred as {\it the Hitchin's hamiltonian system}(\cite{Hitchin1987}):
\[
 T^{*}{\rm Bun}_{G,X}^{D-fl}\stackrel{h}\longrightarrow H^0(X,\,\Omega^{\otimes 2}(D)).
\]
\begin{remark}
By the Riemann-Roch theorem the dimension of $H^0(X,\,\Omega^{\otimes 2}(D))$ is $3(g-1)+N$, which is equal to one of ${\rm Bun}_{G,X}^{D-fl}$.
\end{remark}
Note that $T^{*}{\rm Bun}_{G,X}^{D-fl}$ possesses the canonical
symplectic form. {\bf Remark 2.2} makes us expect the Hitchin's
hamiltonian system should be an algebraically complete integrable
system. We will show this is true.
\par\vspace{5pt}
The inverse image $h^{-1}(0)$ will be referred as {\it the Laumon's
glogal nilpotent variety}(\cite{Laumon}) and will be denoted by $Nilp$.
\begin{prop}
$Nilp$ is a Lagrangian subvariety of $T^{*}{\rm Bun}_{G,X}^{D-fl}$.
\end{prop}
In order to prove the proposition, we will follow the argument of Ginzburg(\cite{Ginzburg}).
Let us fix a Borel subgroup $B$ of $G$ and let $P$ be a principal
$B$-bundle on $X$. By the extension of the structure group we may
associate it a principal $G$-bundle $P_G$:
\[
 P_G=P\times_{B}G.
\]
Note that the natural inclusion
\[
 P \hookrightarrow P_G
\]
yields a $D$-flag by
\[
 P|_{z_i}\hookrightarrow P_G|_{z_i}.
\]
Let ${\rm Bun}_{B,X}$ be the modular stack of principal $B$-bundles on
$X$. The argument above shows that we have a morphism
\[
 {\rm Bun}_{B,X}\stackrel{f}\longrightarrow {\rm Bun}_{G,X}^{D-fl}.
\]
Now we will consider a substack
\[
 Y_f=\{((P,\alpha),(Q,\beta))\in T^{*}{\rm Bun}_{B,X}\times T^{*}{\rm Bun}_{G,X}^{D-fl}\,|\,f(P)=Q,\,f^{*}(\beta)=0\},
\]
and let $p(Y_f)$ be its image under the natural projection
\[
 T^{*}{\rm Bun}_{B,X}\times T^{*}{\rm
 Bun}_{G,X}^{D-fl}\stackrel{p}\longrightarrow T^{*}{\rm
 Bun}_{G,X}^{D-fl}.
\]
\begin{lm}
\[
 p(Y_f)=Nilp.
\]
\end{lm}
{\bf Proof.} Let ${\frak B}$ be the Borel subalgebra corresponding to $B$ and let ${\frak N}$ be its nilpotent ideal.
Then it defines a normal subgroup $N$ of $B$. Let
\[
 Y_{f}\stackrel{\pi}\longrightarrow {\rm Bun}_{B,X}
\]
be a morphism which is defined to be
\[
 \pi((P,\alpha),(Q,\beta))=P.
\]
Then the inverse image $\pi^{-1}(P)$ of $P\in {\rm Bun}_{B,X}$ is
isomorphic to 
\[
 {\rm Ker}[T^{*}_{f(P)}{\rm
 Bun}_{G,X}^{D-fl}\stackrel{f^{*}}\longrightarrow T^{*}_{P}{\rm Bun}_{B,X}].
\]
We will abbreviate 
\[
 H^{i}(\cdot) = H^{i}(X,\,\cdot),
\]
and 
\[
 \cdot\otimes \Omega_X=\cdot\otimes_{{\cal O}_X} \Omega_X. 
\]
Let us consider the diagram:
\[
\begin{array}{ccccccccc}
  &   &  &   &  0 &  & 0 & \\
  &     &  &  & \downarrow & & \downarrow & & \\
  &   &  &   & T^{*}_{f(P)}{\rm
 Bun}_{G,X}^{D-fl} & \stackrel{f^{*}}\to & T^{*}_{P}{\rm Bun}_{B,X} &  & \\
  &     &  &  & \downarrow & & \downarrow & & \\
0 & \to & H^{0}({\rm ad}_{P}({\frak N})(D)\otimes\Omega_X) & \to & H^{0}({\rm ad}_{P}({\frak G})(D)\otimes\Omega_X)  & \to & H^{0}({\rm ad}_{P}({\frak B})(D)\otimes\Omega_X) & \to &0\\
  &     &  &  & \downarrow & & \downarrow & & \\
 &  &  &  & {\frak B}_{D} & \stackrel{id}= & {\frak B}_{D} &  &
\end{array}  
\]
The rightdown vertical arrow is induced by the residue map.
By the standard deformation theory, we have identified
\[
 T^{*}_{P}{\rm Bun}_{B,X}
\]
with 
\[
 H^{0}({\rm ad}_{P}({\frak B})\otimes\Omega_X).
\]
and the right vertical sequence is exact. {\bf Proposition 2.1} shows that
the middle vertical sequence is exact. Finally the middle horizontal sequence is also exact, which
is derived from the sheaf exact sequence
\[
 0 \to {\rm ad}_{P}({\frak N})(D)\otimes\Omega_X \to {\rm
 ad}_{P}({\frak G})(D)\otimes\Omega_X \to {\rm ad}_{P}({\frak
 B})(D)\otimes\Omega_X \to 0, 
\]
and by the fact
\[
 H^{1}({\rm ad}_{P}({\frak N})(D)\otimes\Omega_X)=0.
\]
The diagram chasing shows that
\[
 T^{*}_{f(P)}{\rm
 Bun}_{G,X}^{D-fl}  \stackrel{f^{*}}\longrightarrow  T^{*}_{P}{\rm Bun}_{B,X}
\]
is surjective and that $\pi^{-1}(P)$ is isomorphic to 
\[
 H^{0}({\rm ad}_{P}({\frak N})(D)\otimes\Omega_X).
\]
Hence $p(Y_f)$ is a stack of pairs $(Q,\,\eta)$ of principal $G$-bundles
$Q$ which is induced from a principal $B$-bundle $P$ and $\eta\in H^{0}({\rm
ad}_{P}({\frak N})(D)\otimes\Omega_X)$, which is nothing but $Nilp$.
\begin{flushright}
$\triangle$
\end{flushright}
{\bf Proof of Proposition 2.2.} It is known by Ginzburg
(cf. \cite{Ginzburg} {\bf Lemma 6.5.4}) that $p(Y_f)$ is isotropic. Thus
{\bf Lemma 2.1} shows the dimension of $Nilp$ is less than or equal to
$3(g-1)+N$. But {\bf Remark 2.2} implies 
\[
 \dim Nilp \geq 3(g-1)+N.
\]
Therefore we have obtained the desired result.
\begin{flushright}
$\triangle$
\end{flushright}
The argument of (2.10.1) of \cite{Beilinson-Drinfeld} implies the
following theorem.
\begin{thm}
The Hichtin's hamiltonian system
\[
 T^{*}{\rm Bun}_{G,X}^{D-fl}\stackrel{h}\longrightarrow H^0(X,\,\Omega^{\otimes 2}(D))
\]
is surjective and flat. In particular the dimension of fibres are $3(g-1)+N$.
\end{thm}

\subsection{The Hitchin's correspondence}
We put
\[
 {\cal L}=\Omega_X(D)
\]
and let 
\[
 L\stackrel{\pi}\longrightarrow X
\]
be the ${\mathbb A}^1$-fibration associated to ${\cal L}$. Then the line bundle
$\pi^{*}{\cal L}$ possesses the tautological section $x$. 
An element $q$ of $H^0(X,\,{\cal L}^{\otimes 2})$ determines {\it a
spectral curve} $\Sigma_q$,
which is a closed subscheme of $L$ defined by the equation:
\[
 x^2+\pi^{*}q=0.
\]
This is a double covering of $X$ and its genus $g(\Sigma_q)$ is computed
by the formula of \cite{BNR} {\bf Remark 3.2}:
\begin{equation}
 g(\Sigma_q)=4g-3+N. 
\end{equation}

Let ${\rm Prym}(\Sigma_q\slash X)$ be the {\it Prym variety} of the
double covering $\Sigma_q\stackrel{\pi}\to X$:
\[
 {\rm Prym}(\Sigma_q\slash X)={\rm
 Ker}[{\rm Jac}(\Sigma_q)\stackrel{\pi_*}\longrightarrow {\rm Jac}(X)].
\]
Here ${\rm Jac}$ denotes the Jacobian variety.
The formula (3) shows the dimension of ${\rm Prym}(\Sigma_q\slash X)$ is
equal to $3(g-1)+N$. Also it has the following remarkable property(\cite{BNR}\cite{Hitchin1987}).
\begin{fact}
Suppoose $\Sigma_q$ is smooth. Then there is a bijection between
\[
 {\rm Prym}(\Sigma_q\slash X)
\]
and isomorphism classes of pairs $(E,\,\varphi)$ where $E$ is a vector
 bundle of rank $2$ over $X$ with the trivial determinant and $\varphi$ is a element of
 $H^0(X,\,End_{{\cal O}_X}(E)\otimes{\cal L})$ such that 
\[
 \det (t-\varphi)=t^2+q.
\]									
\end{fact}
Such a bijection will be referred as {\it the Hitchin's correspondence}. We will
regard $H^0(X,\,\Omega_X^{\otimes 2}(D))$ as a subset of $H^0(X,\,{\cal
L}^{\otimes 2})$. Let us take $q\in H^0(X,\,\Omega_X^{\otimes 2}(D))$
so that $\Sigma_q$ is smooth and that the value of $t_i\cdot q$ at each $z_i$ is non-zero.
\begin{prop}
There is a bijection between $h^{-1}(q)$ and ${\rm Prym}(\Sigma_q\slash X)$.
\end{prop}
{\bf Proof.} $h^{-1}(q)$ consists of isomorphism classes of pairs
\[
 ((P,\,\{B_i\}_i),\,\alpha)
\]
where $(P,\,\{B_i\}_i)$ is a principal $G$-bundle on $X$ with a $D$-flag
and 
\[
\alpha\in {\rm Ker}[H^{0}(X,\,{\rm ad}_{P}({\frak G})\otimes{\cal L})\stackrel{\rm Res}\longrightarrow {\frak B}_D].
\]
We claim that $A\in H^{0}(X,\,{\rm ad}_{P}({\frak G})\otimes{\cal L})$ satisfying
\[
 \det A=q
\]
will automatically determine Borel subgroups $\{B_i\}_i$ so that $A$ is
contained in
\begin{equation}
 {\rm Ker}[H^{0}(X,\,{\rm ad}_{P}({\frak G})\otimes{\cal L})\stackrel{\rm Res}\longrightarrow {\frak B}_D].
\end{equation}

In fact, let
\[
 A=
\left(
\begin{array}{cc}
a_i(t_i) & b_i(t_i)\\
c_i(t_i) & -a_i(t_i)
\end{array}
\right)
\frac{dt_i}{t_i}, \quad a_i(t_i), b_i(t_i), c_i(t_i)\in {\mathbb C}[[t_i]], 
\]
be the Taylor expansion of $A$ at $z_i$ and we put
\[
 A_{i}=
\left(
\begin{array}{cc}
a_i(0) & b_i(0)\\
c_i(0) & -a_i(0)
\end{array}
\right).
\]
Since $\det A$ has at most
single pole at $z_i$, $\det A_i$ vanishes. The assumption of residue of
$q$ implies $A_i$ is
a non-zero nilpotent matrix. Now using the Killing form a Borel subalgebra ${\frak B}_i$
is defined to be the orthogonal complement of $A_i$ and let $B_i$ be the
corresponding Borel subgroup. Now it is easy to see that $A$ is
contained in (4). Thus we have proved there is one to one correspondence
between $h^{-1}(q)$ and isomorphism classes of pairs of a principal
$G$-bundle $P$ on $X$ and  $A\in H^{0}(X,\,{\rm ad}_{P}({\frak G})\otimes{\cal L})$ satisfying
\[
 \det A=q.
\] 
Now the natural correspondence between isomorphism classes of principal
$G$-bundles on $X$ and of rank 2 vector bundles on $X$ with the trivial
determinant and {\bf Fact 2.1} imply the desired result. 
\begin{flushright}
$\triangle$
\end{flushright}
{\bf Theorem 2.1} and {\bf Proposition 2.3} will show the following theorem.
\begin{thm}
The Hichtin's hamiltonian system
\[
 T^{*}{\rm Bun}_{G,X}^{D-fl}\stackrel{h}\longrightarrow H^0(X,\,\Omega^{\otimes 2}(D))
\]
is surjective and flat. Moreover there is a Zariski open dense subset
 $U$ of $H^0(X,\,\Omega^{\otimes 2}(D))$ over which $h$ is an
 abelian fibration.
\end{thm}
\begin{cor}
The Hichtin's hamiltonian system induces an isomorphism of commutative  ${\mathbb C}$-algebras:
\[
 \Gamma(H^0(X,\,\Omega^{\otimes 2}(D)),\,{\cal O}) \stackrel{h^{*}}\simeq \Gamma(T^{*}{\rm Bun}_{G,X}^{D-fl},\,{\cal O}).
\]
\end{cor}

\section{A uniformization of a modular stack}
Throughout the paper, we will use the following notations. {\it The loop
group} $LG$ and {\it the regular loop group} $L_+G$ are defined to be
\[
 LG=G({\mathbb C}((t))), \quad L_+G=G({\mathbb C}[[t]]), 
\]
respectively. Let ${\rm ev}_0$ be the evaluation map:
\[
 L_+G \stackrel{{\rm ev}_0} \longrightarrow G, \quad {\rm ev}_0(A)=A(0),
\]
and we set
\[
 L_{0}G={\rm Ker}[L_+G \stackrel{{\rm ev}_0} \longrightarrow G].
\]
Let $B_+$ be the upper Borel subgroup of $G$ and we put
\[
 \hat{B}_+={\rm ev}_0^{-1}(B_+).
\]
Also we will consider the corresponding Lie algebras:
\begin{itemize}
\item $L{\frak G}={\frak G}\otimes_{\mathbb C} {\mathbb C}((t))$ : the loop algebra.
\item $L_+{\frak G}={\frak G}\otimes_{\mathbb C} {\mathbb C}[[t]]$ : the regular
      loop algebra.
\item $L_0{\frak G}={\frak G}\otimes_{\mathbb C} t{\mathbb C}[[t]]$.
\item $\hat{\frak B}_+={\frak B}_+\oplus L_0{\frak G}$.
\end{itemize}

\subsection{The vacuum module of the critical level}
We will consider a central extension
\begin{equation}
0 \to {\mathbb C}K \to \hat{\frak G}_N \to (L{\frak G})^{\oplus N} \to 0
\end{equation}
which is defined in the following way. Let 
\[
 A=(A_1\otimes f_1,\cdots, A_N\otimes f_N), \quad  B=(B_1\otimes g_1,\cdots, B_N\otimes g_N),
\]
where
\[
 A_i,\,B_j \in {\frak G},\quad f_i,\,g_j \in {\mathbb C}((t))
\]
be elements of $(L{\frak G})^{\oplus N}$. Then their commutator is defined to
be
\[
 [A,\,B]=([A_1,B_1]\otimes f_1g_1,\cdots,[A_N,B_N]\otimes
 f_Ng_N)-\{\sum_{i=1}^N{\rm Tr}(A_iB_i)\cdot {\rm
 Res}_{t=0}(f_i\cdot dg_i)\}\cdot K.
\]
When $N=1$ we set 
\[
 \hat{\frak G}=\hat{\frak G}_1.
\]
Obviously (5) splits over $(L_+{\frak G})^{\oplus N}$. We put 
\[
 Y=X\setminus \{z_1,\cdots,z_N\},
\]
and 
\[
 {\frak G}(Y)={\frak G}\otimes_{\mathbb C}\Gamma(Y,\,{\cal O}),
\]
which is imbedded in $(L{\frak G})^{\oplus N}$ by the Taylor
expansion. Let us observe that (5) also splits over ${\frak G}(Y)$ by
the residue theorem. Hence we may consider them as Lie subalgebra of
$\hat{\frak G}_N$.
\par\vspace{5pt}
Let $U(\hat{\frak G}_N)$ be the universal envelopping algebra of
$\hat{\frak G}_N$ and its {\it reduction at the critical level} is
defined to be
\[
 U_{-2}(\hat{\frak G}_N)=U(\hat{\frak G}_N)\slash(K+2).
\]
Let ${\mathbb C}_{-2}$ be a ${\mathbb C}K\oplus (L_0{\frak G})^{\oplus
N}$-module which is isomorphic to ${\mathbb C}$ as an abstract vector space
whose action is defined to be
\[
 (L_0{\frak G})^{\oplus N}\cdot {\mathbb C}_{-2}=0
\] 
and
\[
 K=-2.
\]
Then we define {\it the vacuum module} $Vac_{-2}^{0}$ as
\[
Vac_{-2}^{0}=U_{-2}(\hat{\frak G}_N)\otimes_{U({\mathbb C}K\oplus (L_0{\frak
G})^{\oplus N})}{\mathbb C}_{-2},
\]
and let $v_{-2}$ be the the highest weight vector.
\subsection{A uniformization of a modular stack}

Let ${\rm Bun}_{G,X}^{D-fr}$ be the modular stack over $Aff_{\slash{\mathbb
C}, et}$ which classifies isomorphism classes of pairs
\[
(P,\,\{\alpha_i\}), 
\]
where $P$ is a principal $G$-bundle on $X$ and $\alpha_i$ is its
trivialization $z_i$:
\[
 G\stackrel{\alpha_i}\simeq P|_{z_i}.
\]
(Such an $N$-tuple of trivializations will be called as {\it a $D$-framing of} $P$.)
This becomes a principal $(B_+)^{N}$-bundle over ${\rm
Bun}_{G,X}^{D-fl}$ by the morphism
\[
 {\rm Bun}_{G,X}^{D-fr}\stackrel{\pi}\longrightarrow {\rm Bun}_{G,X}^{D-fl},
\]
where
\[
 \pi((P,\,\{\alpha_i\}))=(P,\,\{\alpha_i(B_+)\}).
\]
By the Taylor expansion 
\[
 G(Y)=G(\Gamma(Y,\,{\cal O}))
\]
may be considered as a subgroup of $(LG)^{N}$ and we set
\[
 Z=(LG)^{N}\slash G(Y).
\]
Then each of ${\rm Bun}_{G,X}^{D-fr}$, ${\rm Bun}_{G,X}^{D-fl}$ and ${\rm
Bun}_{G,X}$ has the following uniformization as a stack in terms of
$Z$(\cite{Beilinson-Drinfeld} \cite{LS} \cite{Sorger}):
\[
 {\rm Bun}_{G,X}^{D-fr}\simeq (L_0G)^{N}\backslash Z,\quad {\rm
 Bun}_{G,X}^{D-fl}\simeq (\hat{B}_+)^{N}\backslash Z,
\]
and 
\[
 {\rm Bun}_{G,X}\simeq (L_+G)^{N}\backslash Z.
\]
Let $\Theta_Z$ be the tangent sheaf of $Z$. Then the left action of
$(LG)^N$ on $Z$ induces a surjection
\[
 (L{\frak G})^{\oplus N}\otimes_{\mathbb C}{\cal O}_Z \longrightarrow \Theta_Z.
\]
Combining with this the projection
\[
 \hat{\frak G}_N \longrightarrow (L{\frak G})^{\oplus N},
\]
we obtain a surjective map
\[
 \hat{\frak G}_N\otimes_{\mathbb C}{\cal O}_Z
 \stackrel{\alpha}\longrightarrow \Theta_Z,
\]
which will referred as {\it the anchor
map}(\cite{Beilinson-Bernstein1993}, \cite{Frenkel-BenZvi} {\bf
Appendix 3}). 
We will consider a ${\cal O}_Z$-module
\[
 {\cal J}={\hat {\frak G}}_N\otimes_{\mathbb C}{\cal O}_Z\slash {\rm Ker}\,\alpha,
\]
which is a central extension of $\Theta_Z$:
\[
0 \to {\cal O}_Z \to {\cal J} \to \Theta_Z \to 0.
\]
Now we set
\[
 {\cal D}^{\prime}_Z=U({\cal J})\slash(K+2).
\]
This is a sheaf of ${\mathbb C}$-algebra on $Z$ and the anchor map defines
an algebra homomorphism
\begin{equation}
 U_{-2}(\hat{\frak G}_N) \stackrel{U(\alpha)}\longrightarrow {\cal D}^{\prime}_Z,
\end{equation}
which is compatible with the action of $(L_+{\frak G})^{\oplus N}$.
\subsection{A sheaf of twisted differential operators on ${\rm
  Bun}_{G,X}^{D-fl}$}
Let us fix a theta characteristic $\kappa$ of $X$. Namely we will fix a line
bundle on $X$ so that $\kappa^{\otimes 2}\simeq \Omega_X$. Then it will
determine a half canonical $\omega^{\otimes \frac{1}{2}}$ of ${\rm
Bun}_{G,X}$(\cite{Sorger}).
\par\vspace{5pt}
Let
\[
 Z\stackrel{\tau}\longrightarrow {\rm Bun}_{G,X}
\]
be the projection. Then it is
known that $(\tau_{*}{\cal
D}_{Z}^{\prime})^{(L_+G)^N}$ is isomorphic to the sheaf of differential operators
${\cal D}^{\prime}$ on $\omega^{\otimes \frac{1}{2}}$
(\cite{Beilinson-Drinfeld}\cite{Frenkel1995}).

For 
\[
 \lambda=(\lambda_1,\cdots,\lambda_N)\in {\mathbb Z}^N,
\]
let ${\cal O}(\lambda)$ be the line bundle on $({\mathbb P}^1)^N$ which is
defined to be
\[
 {\cal O}(\lambda)=p_1^{*}{\cal O}_{{\mathbb P}^1}(\lambda_1)\otimes \cdots \otimes
 p_N^{*}{\cal O}_{{\mathbb P}^1}(\lambda_N),
\]
where $p_i$ is the projection to the $i$-th factor. Let $S$ be an affine scheme
defined over ${\mathbb C}$ and let
\[
 X\times S \stackrel{p_S}\longrightarrow S
\]
be the projection. Then a principal $G$-bundle $P$ on $X\times S$ induces a
Cartesian diagram
\[
\begin{array}{ccc}
({\mathbb P}^1)^N\times S   & \longrightarrow & {\rm Bun}_{G,X}^{D-fl}\\
\downarrow & & \downarrow\\
S & \stackrel{P}\longrightarrow & {\rm Bun}_{G,X}.
\end{array}
\]
By definition, the restriction $\omega|_S$ of the canonical sheaf
$\omega$ of ${\rm Bun}_{G,X}$ to $S$ is given by
\[
 \omega|_S={\rm det}(Rp_{S*}{\rm ad}_{P}(\frak G))^{\otimes (-1)},
\]
and in particular we have
\[
 \omega^{\otimes \frac{1}{2}}|_S={\rm det}(Rp_{S*}{\rm ad}_{P}(\frak G))^{\otimes (-\frac{1}{2})}.
\]
(Note that in order to define the right hand side, it is necessary to
choose a theta characteristic.)
\par\vspace{5pt}
Now we define a line bundle ${\cal L}_{\lambda}$ on ${\rm
Bun}_{G,X}^{D-fl}$ so that its restriction to $({\mathbb P}^1)^N\times S$ is
equal to 
\[
 p^{*}{\cal O}(\lambda)\otimes q^{*}(\omega^{\frac{1}{2}}|_S),
\]
where $p$ (resp. $q$) be the projection onto $({\mathbb P}^1)^N$
(resp. $S$).
The sheaf of differential operators on ${\cal L}_{\lambda}$ will be
denoted by ${\cal D}^{\prime}_{D-fl,\lambda}$ and let
$D^{\prime}_{D-fl,\lambda}$ be its global sections. As before it has the
following description in terms of the uniformization.
\par\vspace{5pt}
Let
\[
 Z\stackrel{\tau_{D-fr}}\longrightarrow {\rm Bun}^{D-fr}_{G,X}
\]
be the projection and we put (\cite{Frenkel-BenZvi} {\bf 16.2}):
\[
 {\cal D}^{\prime}_{D-fr}=(\tau_{D-fr*}{\cal D}_Z^{\prime})^{(L_0G)^N}.
\]
This is a sheaf of differential operators on ${\rm Bun}_{G,X}^{D-fr}$
and let $D^{\prime}_{D-fr}$ be the space of its global sections. Since
$U({\frak G})^{\otimes N}$ is a subalgebra of $U_{-2}(\hat{\frak G})$,
(6) induces a homomorphism:
\[
 U({\frak G})^{\otimes N}\longrightarrow {\cal D}_Z^{\prime}.
\]
Note that $U({\frak G})^{\otimes N}$ is invariant under the action of
$(L_0G)^N$ and we have a homomorphism
\begin{equation}
U({\frak G})^{\otimes N}\longrightarrow D^{\prime}_{D-fr}.
\end{equation}
For an integer $k$, let $M_k$ be the Verma module of ${\frak G}$ of the
highest weight $k$ with the highest weight vector $\mu_k$. We will consider a $U({\frak G})^{\otimes N}$-module
\[
 M_{\lambda}=M_{\lambda_1}\otimes \cdots \otimes M_{\lambda_N},
\]
and let $\mu_{\lambda}$ be its highest weight vector:
\[
 \mu_{\lambda}=\mu_{\lambda_1}\otimes \cdots \otimes \mu_{\lambda_N}.
\]
Using the map (7), one has the sheaf of twisted differential operators
on ${\rm Bun}_{G,X}^{D-fl}$
\[
 [\pi_{*}({\cal D}^{\prime}_{D-fr}\otimes_{U({\frak G})^{\otimes N}}M_{\lambda})]^{B^N},
\]
which is easily seen to be isomorphic to ${\cal
D}^{\prime}_{D-fl,\lambda}$(cf. \cite{Frenkel1995} {\bf 5.2}).
\section{The Beilinson-Drinfeld homomorphism}
\subsection{The Feigin-Frenkel isomorphism}
We put
\[
 I=(L_0{\frak G})^{\oplus N},
\]
which will be considered as a Lie subalgebra of $\hat{{\frak G}}_N$.
Let $\tilde{{\frak D}}_{{\frak G},N}$ be the left normalizer of
$U_{-2}(\hat{{\frak G}}_N)\cdot I$ in $U_{-2}(\hat{{\frak G}}_N)$ and we
set
\[
 {\frak D}_{{\frak G},N}=\tilde{{\frak D}}_{{\frak G},N}\slash U_{-2}(\hat{{\frak G}}_N)\cdot I.
\]
Note that $G^{N}$ preserves both of $U_{-2}(\hat{{\frak
G}}_N)\cdot I$ and $\tilde{{\frak D}}_{{\frak G},N}$ and therefore
it acts on ${\frak D}_{{\frak G},N}$.
The map (6) induces a homomorphism(cf. \cite{Beilinson-Drinfeld}{\bf 1.1.3})
\begin{equation}
 {\frak D}_{{\frak G},N}\longrightarrow D^{\prime}_{D-fr}, 
\end{equation}
which is compatible with the $G^{N}$ action.
Putting 
\[
 {\rm deg}\,X=1
\]
for any $X\in (L{\frak G})^{\oplus N}$, $U_{-2}(\hat{{\frak G}}_N)$
becomes a filtered algebra and ${\frak D}_{{\frak G},N}$ has the induced
filtration. On the other hand $D^{\prime}_{D-fr}$ possesses a filtration
by the order of a differential operator and (8) preserves these
filtrations. 
\par\vspace{5pt}
Let $\Sigma(Vac_{-2}^0)$ be the space of singular vectors:
\[
 \Sigma(Vac_{-2}^0)=\{ x\in Vac_{-2}^0 \,| \,\gamma\cdot x=0
 \quad\mbox{for}\quad \gamma\in I\}.
\]
Then the map
\[
 U_{-2}(\hat{{\frak G}}_N)\stackrel{\epsilon}\longrightarrow Vac_{-2}^0,
 \quad \epsilon(X)=X\cdot v_{-2},
\]
induces an isomorphism of vector spaces:
\[
 {\frak D}_{{\frak G},N}\stackrel{\epsilon^{\prime}}\simeq \Sigma(Vac_{-2}^0).
\]
On the other hand, let $Z_{-2}({\hat{\frak G}})$ be the center of the
local completion of $U_{-2}(\hat{{\frak G}}_N)$(\cite{Frenkel02}). Then
$Z_{-2}({\hat{\frak G}})^{\otimes N}$ acts on $Vac_{-2}^0$ and since
$Z_{-2}({\hat{\frak G}})$ commutes with $U_{-2}(\hat{{\frak G}}_N)$ we
have a map
\[
 Z_{-2}({\hat{\frak G}})^{\otimes N}\stackrel{\rho}\longrightarrow \Sigma(Vac_{-2}^0),
\]
which is defined to be
\[
 \rho(Z)=Z\cdot v_{-2},\quad Z\in  Z_{-2}({\hat{\frak G}})^{\otimes N}.
\]
Thus we have a linear map:
\begin{equation}
 Z_{-2}({\hat{\frak G}})^{\otimes N}\stackrel{{\epsilon}^{\prime
 -1}\circ\rho}\longrightarrow {\frak D}_{{\frak G},N},
\end{equation}
which is checked to be a homomorphism of ${\mathbb C}$-algebras.
Note that the map
(9) is compatible with the action of $G^{N}$. In particular since $Z_{-2}({\hat{\frak
G}})$ is the center, its image is contained in
\[
 {\frak D}_{{\frak G},N}^{B^N}=\{x\in{\frak D}_{{\frak G},N}\,|\,b\cdot
 x=x\quad \mbox{for}\quad b\in B^N \}.
\]
By \cite{Feigin-Frenkel} (see also \cite{Frenkel02} $\S12$), it is known
that
$Z_{-2}({\hat{\frak G}})$ is isomorphic to a topological completion of a
polynomial algebra ${\mathbb C}[\{S_m\}_{m\in{\mathbb Z}}]$ generated by
the Sugawara operators $\{S_m\}_{m\in{\mathbb Z}}$. Thus we have obtained a
homomorphism 
\begin{equation}
 {\mathbb C}[\{S_m\}_{m\in{\mathbb Z}}]^{\otimes N}\longrightarrow {\frak D}_{{\frak G},N}^{B^N}.
\end{equation}
Since each $S_{m}$ is an infinite sum of normally ordered quadratics of
$L{\frak G}$(\cite{Frenkel-BenZvi} {\bf 2.5.10}), if we put
\[
 {\rm deg}\,S_m=2,
\]
and give a filtration ${\mathbb C}[\{S_m\}_{m\in{\mathbb Z}}]$ by the
degree, (10) becomes a filtered algebra homomorphism.
\par\vspace{5pt}
We will use the following notations:
\begin{itemize}
\item $D={\rm Spec}\,{\mathbb C}[[t]]$ : a formal disc,
\item $D^{\times}={\rm Spec}\,{\mathbb C}((t))$ : a punctured formal disc.
\end{itemize}

Let $c_m$ be a linear functional on $H^{0}(D^{\times},\,\Omega^{\otimes
2})$ which is defined to be
\[
 c_m(q)=q_m,\quad \mbox{for}\quad q=\sum_{n\in{\mathbb
 Z}}q_nt^{-n-2}(dt)^{\otimes 2}\in H^{0}(D^{\times},\,\Omega^{\otimes
2}).
\]
Then ${\mathbb C}[\{S_m\}_{m\in{\mathbb Z}}]$ is isomorphic to
\[
 \Gamma(H^{0}(D^{\times},\,\Omega^{\otimes
2}),\,{\cal O})\simeq {\mathbb C}[\{c_m\}_{m\in{\mathbb Z}}],
\]
by a map
\[
 \phi(S_m)=c_m.
\]
Thus (10) turns out to be a homomorphism:
\[
\Gamma(H^{0}(D^{\times},\,\Omega^{\otimes
2}),\,{\cal O})^{\otimes N}\stackrel{\psi}\longrightarrow {\frak D}_{{\frak G},N}^{B^N}.  
\]

For $S\in Z_{-2}({\hat{\frak G}})$, we put
\[
 S^{(i)}=1\otimes \cdots \otimes S \otimes \cdots \otimes 1\in
 Z_{-2}({\hat{\frak G}})^{\otimes N},
\]
where $S$ appears in the $i$-th place. Then by the reason of degree, we
will find
\[
 S_m^{(i)}\cdot v_{-2}=0,
\]
for any $1\leq i \leq N$ and $m > 0$. Thus the homomorphism above factors through 
\begin{equation}
 \Gamma(H^{0}(D,\,\Omega^{\otimes
2}(2{\mathbb O})),\,{\cal O})^{\otimes N}\stackrel{\psi_+}\longrightarrow {\frak D}_{{\frak G},N}^{B^N}.
\end{equation}

Here, in general for an integer $k$, $H^{0}(D,\,\Omega^{\otimes
2}(k{\mathbb O}))$ will stand for the space of quadratic differentials on the formal
disc which has a pole at the origin ${\mathbb O}$ at most $k$-th order. Note
that $\Gamma(H^{0}(D,\,\Omega^{\otimes
2}(2{\mathbb O})),\,{\cal O})$ is isomorphic to ${\mathbb C}[\{c_m\}_{m\leq
0}]$. 
Putting
\[
 {\rm deg}\,c_m=2,
\]
for any $m$, we will introduce a filtration on $\Gamma(H^{0}(D,\,\Omega^{\otimes
2}(2{\mathbb O})),\,{\cal O})\simeq {\mathbb C}[\{c_m\}_{m\leq
0}]$.
Then (11) preserves the filtrations and combining it with (8),
we will obtain a filtered algebra homomorphism
\begin{equation}
  \Gamma(H^{0}(D,\,\Omega^{\otimes
2}(2{\mathbb O})),\,{\cal O})^{\otimes N}\stackrel{\beta}\longrightarrow (D^{\prime}_{D-fr})^{B^N}.
\end{equation}

\subsection{The Beilinson-Drinfeld homomorphism}
By definition of the Verma module, there is a surjection as $U({\frak
G})^{\otimes N}$-modules:
\begin{equation}
 U({\frak G})^{\otimes N} \stackrel{\rho_{\lambda}}\longrightarrow M_{\lambda},
\end{equation}
which is defined to be
\[
 \rho_{\lambda}(X_{1}\otimes \cdots \otimes X_N)=(X_{1}\otimes \cdots
 \otimes X_N)\cdot \mu_{\lambda}.
\]
The map (13) induces a homomorphism
\[
 {\cal D}^{\prime}_{D-fr}\simeq {\cal D}^{\prime}_{D-fr}\otimes_{U({\frak
 G})^{\otimes N}}U({\frak G})^{\otimes N} \longrightarrow {\cal
 D}^{\prime}_{D-fr}\otimes_{U({\frak G})^{\otimes N}}M_{\lambda},
\]
which implies (cf. {\bf 3.3})
\[
 (\pi_*{\cal D}^{\prime}_{D-fr})^{B^N} \longrightarrow [\pi_*({\cal
 D}^{\prime}_{D-fr}\otimes_{U({\frak G})^{\otimes
 N}}M_{\lambda})]^{B^N}\simeq {\cal D}^{\prime}_{D-fl, \lambda}.
\]
Taking global sections on ${\rm Bun}_{G,X}^{D-fl}$, we have a
homomorphism 
\begin{equation}
 (D^{\prime}_{D-fr})^{B^N} \longrightarrow D^{\prime}_{D-fl, \lambda}.
\end{equation}
Note that, by the construction, this preserves the
filtrations. Combining (14) with (12), we will obtain the following
proposition.
\begin{prop}
There is a homomorphism of algebras which preserves the filtrations
\[
 \Gamma(H^0(D,\,\Omega^{\otimes 2}(2{\mathbb O})),\,{\cal O})^{\otimes N} \stackrel{\hat{\beta}_{\lambda}}\longrightarrow D^{\prime}_{D-fl, \lambda}.
\]
\end{prop}
For $1 \leq i \leq N$, we set
\[
 c_{m}^{(i)}=1\otimes \cdots \otimes 1 \otimes c_{m} \otimes 1\otimes
 \cdots \otimes 1 \in \Gamma(H^0(D,\,\Omega^{\otimes 2}(2{\mathbb O})),\,{\cal O})^{\otimes N},
\]
where $c_{m}$ appears in the $i$-th place. Then it is easily checked
(\cite{Frenkel1995} {\bf 3.6} Example) that
\[
 \hat{\beta}_{\lambda}(c_{0}^{(i)})=\frac{\lambda_i(\lambda_i+2)}{4}.
\]
Thus {\bf Proposition 4.1} implies the following theorem.
\begin{thm}
There is a homomorphism of algebras which preserves the filtrations
\[
 \Gamma(H^0(D,\,\Omega^{\otimes 2}({\mathbb O})),\,{\cal O})^{\otimes N} \stackrel{{\beta}_{\lambda}}\longrightarrow D^{\prime}_{D-fl, \lambda}.
\]
\end{thm}
The homomorphism ${\beta}_{\lambda}$ will be referred as {\it the
Beilinson-Drinfeld homomorphism}.

\section{A quantization of the Hitchin system}
\subsection{A quantization of the Hitchin's hamiltonian system}
Let 
\[
 \{f_1,\cdots,f_{3(g-1)+N}\}
\]
be a basis of $H^0(X,\,\Omega^{\otimes 2}(D))$ and for each $i$ we put
\[
 {\rm deg}\,f_i=2.
\]
Then $\Gamma(H^0(X,\,\Omega^{\otimes 2}(D)),\,{\cal O})$ is isomorphic
to a weighted polynomial ring
\[
 {\mathbb C}[f_1,\cdots,f_{3(g-1)+N}].
\]
Let $D_i$ be the formal disc at $z_i$:
\[
 D_i={\rm Spec}\,{\mathbb C}[[t_i]].
\]
The restriction map
\[
 H^0(X,\,\Omega^{\otimes 2}(D))\stackrel{r}\longrightarrow
 \oplus_{i=1}^N H^0(D_i,\,\Omega^{\otimes 2}({\mathbb O}))
\]
induces a surjective homomorphism
\[
 \otimes_{i=1}^{N}\Gamma(H^0(D_i,\,\Omega^{\otimes 2}({\mathbb O})),\,{\cal
 O}) \stackrel{r^*}\longrightarrow \Gamma(H^0(X,\,\Omega^{\otimes
 2}(D)),\,{\cal O}),
\]
which preserves the degree. As \cite{Beilinson-Drinfeld} {\bf 2.4.3} one
may see that the Beilinson-Drinfeld homomorphism yields
a commutative diagram:
\[
 \begin{array}{ccc}
\otimes_{i=1}^{N}\Gamma(H^0(D_i,\,\Omega^{\otimes 2}({\mathbb O})),\,{\cal
 O}) & \stackrel{{\rm Gr}({\beta_{\lambda}})}\longrightarrow & \Gamma(T^*{\rm Bun}_{G,X}^{D-fl},\,{\cal O})\\
\Vert &  & \uparrow h^{*} \\
\otimes_{i=1}^{N}\Gamma(H^0(D_i,\,\Omega^{\otimes 2}({\mathbb O})),\,{\cal
 O}) & \stackrel{r^*}\longrightarrow & \Gamma(H^0(X,\,\Omega^{\otimes
 2}(D)),\,{\cal O})
\end{array}
\]
Here the right vertical arrow $h^{*}$ is induced by the Hitchin's
hamiltonian system, which is an isomorphism as we have seen at {\bf
Corollary 2.1}. Since $r^{*}$ is surjective, so is ${\rm
Gr}({\beta_{\lambda}})$. An easy induction argument with respect to the
filtration will show that $\beta_{\lambda}$ is surjective. In
particular, we have obtained the following proposition.
\begin{prop}
$D^{\prime}_{D-fl,\lambda}$ is a commutative ring.
\end{prop}
Let $F_{\cdot}$ be the filtration.
For each $f_i$, we choose $\delta_i \in F_2(D^{\prime}_{D-fl,\lambda})$
so that
\[
 [\delta_i]=h^{*}f_i \in {\rm
 Gr}_{2}^{F}(D^{\prime}_{D-fl,\lambda})\simeq \Gamma({\rm
 Bun}_{G,X}^{D-fl},\,{\rm Sym}^2(T^*{\rm
 Bun}_{G,X}^{D-fl})),
\]
where ${\rm Sym}^2(T^*{\rm
 Bun}_{G,X}^{D-fl})$ is the symmetric square of the cotangent bundle of ${\rm
 Bun}_{G,X}^{D-fl}$. This defines a homomorphism
\[
 \Gamma(H^0(X,\,\Omega^{\otimes
 2}(D)),\,{\cal O}) \stackrel{h_{\lambda}}\longrightarrow
 D^{\prime}_{D-fl,\lambda},\quad h_{\lambda}(f_i)=\delta_i.
\]
which preserves the filtrations and that
\[
 {\rm Gr}(h_{\lambda})=h^{*}.
\]
As before, by the induction argument with respect to the filtration, one can see
that $h_{\lambda}$ is surjective.
\begin{lm}
$h_{\lambda}$ is injective.
\end{lm}
{\bf Proof.} If $h_{\lambda}$ were not injective, there is 
\[
 x\neq 0\in \Gamma(H^0(X,\,\Omega^{\otimes
 2}(D)),\,{\cal O})
\]
such that $h_{\lambda}(x)=0$. We take the minimal nonnegative integer
$i$ so that
\[
 x\in F_i\setminus F_{i-1}.
\]
Here $F_i$ denotes $F_{i}(\Gamma(H^0(X,\,\Omega^{\otimes
 2}(D)),\,{\cal O}))$. Let $[x]$ be the element of ${\rm Gr}_{i}^{F}((\Gamma(H^0(X,\,\Omega^{\otimes
 2}(D)),\,{\cal O})))$ determined by $x$. By our choice of $i$, we have
\begin{equation}
 [x]\neq 0.
\end{equation}
Now $h_{\lambda}(x)=0$ implies
\[
 h^*([x])={\rm Gr}(h_{\lambda})([x])=0.
\]
Since $h^{*}$ is an isomorphism, $[x]$ must be $0$. But this contradicts to (15). 
\begin{flushright}
$\triangle$
\end{flushright}
Thus we have proved the following theorem.
\begin{thm}
There is an isomorphism of filtered algebras:
\[
 \Gamma(H^0(X,\,\Omega^{\otimes
 2}(D)),\,{\cal O}) \stackrel{h_{\lambda}}\simeq D^{\prime}_{D-fl,\lambda}, 
\]
such that
\[
 {\rm Gr}(h_{\lambda})=h^*.
\]
\end{thm}
Let 
\[
 \Omega^{\otimes 2}(2D)\stackrel{R_{-2}}\longrightarrow {\mathbb C}^N
\]
be a map which is defined to be
\[
 R_{-2}(\omega)=({\rm Res}_{t_1=0}(t_1\cdot \omega),\cdots,{\rm
 Res}_{t_N=0}(t_N\cdot \omega)),\quad \omega \in \Omega^{\otimes 2}(2D).
\]
Then we have an exact sequence
\begin{equation}
 0 \to \Omega^{\otimes 2}(D) \to \Omega^{\otimes 2}(2D) \stackrel{R_{-2}}\to {\mathbb C}^N
 \to 0.
\end{equation}
Note that the assumption for $\Omega(D)$ implies
\[
 H^{1}(X,\,\Omega^{\otimes 2}(D))=0.
\]
Thus (16) induces an exact sequence
\[
 0 \to H^0(X,\,\Omega^{\otimes 2}(D)) \to H^0(X,\,\Omega^{\otimes
 2}(2D)) \stackrel{R_{-2}}\to {\mathbb C}^N \to 0.
\]
Now we put
\[
 \Delta(\lambda_i)=\frac{\lambda_{i}(\lambda_{i}+2)}{4},
\]
and
\[
 \Delta(\lambda)=(\Delta(\lambda_1),\cdots, \Delta(\lambda_N))\in {\mathbb C}^N.
\]
Let us consider an affine variety
\[
 H_{\Delta(\lambda)}=(R_{-2})^{-1}(\Delta(\lambda)).
\]
This is a homogeneous space of $H^0(X,\,\Omega^{\otimes 2}(D))$ and the
translation by $\Delta(\lambda)$ gives an isomorphism between $\Gamma(H_{\Delta(\lambda)},\,{\cal O})$ and
$\Gamma(H^0(X,\,\Omega^{\otimes 2}(D)),\,{\cal O})$ as filtered
algebras. Now {\bf Theorem 5.1} implies
\begin{thm}
There is an isomorphism of ${\mathbb C}$-algebras:
\[
  \Gamma(H_{\Delta(\lambda)},\,{\cal O})\stackrel{h_{\lambda}^{+}}\simeq D^{\prime}_{D-fl,\lambda}, 
\]
which preserves the filtrations and that
\[
 {\rm Gr}(h_{\lambda}^{+})=h^*.
\]
\end{thm}
The isomorphism $h_{\lambda}^{+}$ will be referred as {\it a
quantization of the Hitchin's hamiltonian system}.
\subsection{A generalization of the Beilinson-Drinfeld correspondence}
An element $q$ of $H_{\Delta(\lambda)}$ defines a homomorphism of
algebra:
\[
 \Gamma(H_{\Delta(\lambda)},\,{\cal O})\stackrel{f_q}\longrightarrow
 {\mathbb C}.
\]
By {\bf Thorem 5.2}, ${\rm Ker}\, f_q$ may be considered as a kernel of
$D^{\prime}_{D-fl,\lambda}$. We will define a D-module ${\cal M}_q$ on
${\rm Bun}_{G,X}^{D-fl}$ to be
\[
 {\cal M}_q=({\cal D}^{\prime}_{D-fl,\lambda}\slash {\rm Ker}\,
 f_q)\otimes_{{\cal O}}{\cal L}_{\lambda}^{-1}.
\]
One may see that the chracteristic variety(\cite{Borel}) $char({\cal
M}_q)$ of ${\cal M}_q$ coinsides with
$Nilp$(cf. \cite{Beilinson-Drinfeld}{\bf 5.1.2 Proposition}). {\bf
Proposition 2.2} implies the following theorem.
\begin{thm}
 ${\cal M}_q$ is a holonomic D-module.
\end{thm}
The correspondence which assigns the D-module ${\cal M}_q$ to $q\in
H_{\Delta(\lambda)}$ will be referred as {\it a generalized
Beilinson-Drinfeld correspondence}.

\section{A localization functor}
Let
\[
 Z \stackrel{\tau_{D-fl}}\longrightarrow {\rm Bun}_{G,X}^{D-fl}\simeq
 \hat{B}_+^N\backslash Z
\]
be the projection. Using the map (6), one can associate a D-module
$\Delta_{D-fl}(M)$ on ${\rm Bun}_{G,X}^{D-fl}$ to a $U_{-2}(\hat{{\frak G}}_N)$-module $M$ by
\[
 \Delta_{D-fl}(M)=[\tau_{D-fl*}({\cal D}_{Z^{\prime}}\otimes_{U_{-2}(\hat{{\frak G}}_N)}M)]^{\hat{B}_+^N}.
\]
The functor $\Delta_{D-fl}$ will be referred as {\it a localization
functor} (\cite{Frenkel-BenZvi} {\bf 16.2}). Let us remember that ${\frak G}^{\oplus N}$ is a Lie subalgebra of
$\hat{{\frak G}}_N$. Therefore we have a homomorphism
\[
 U({\frak G})^{\otimes N} \longrightarrow U_{-2}(\hat{{\frak G}}_N),
\]
and for $\lambda=(\lambda_1,\cdots,\lambda_N)\in {\mathbb Z}^N$, we may form
$M_{\lambda}$ into $U_{-2}(\hat{{\frak G}}_N)$-module $M_{-2,\lambda}$:
\[
 M_{-2,\lambda}=U_{-2}(\hat{{\frak G}}_N)\otimes_{U({\frak G})^{\otimes N}}M_{\lambda}.
\]
As we have seen in {\bf 3.3}, $\Delta_{D-fl}(M_{-2,\lambda})$ is
isomorphic to ${\cal D}^{\prime}_{D-fl,\lambda}$.
\begin{df}
We will call
\[
 q=(q_1,\cdots,q_N)\in \oplus_{i=1}^NH^0(D_i,\,\Omega^{\otimes 2}(2{\mathbb O}))
\]
is {\rm global} if it is contained in the image of the Taylor expansion
\[
 H^0(X,\,\Omega^{\otimes 2}(2D)) \longrightarrow \oplus_{i=1}^NH^0(D_i,\,\Omega^{\otimes 2}(2{\mathbb O})).
\]
Also it will be called $\lambda$-{\rm admissible} if it has a Taylor
 expansion
\[
 q_i=q_i(t_i)dt^{\otimes 2}
\]
such that
\[
 q_i(t_i)=\frac{\Delta(\lambda_i)}{t_i^2}+\cdots,
\]
for each $i$.
\end{df}
Let 
\[
 q=(q_1,\cdots,q_N)\in \oplus_{i=1}^NH^0(D_i,\,\Omega^{\otimes 2}(2{\mathbb O}))
\]
be $\lambda$-admissible. It defines a homomorphism
\[
 \otimes_{i=1}^N\Gamma(H^0(D_i,\,\Omega^{\otimes 2}(2{\mathbb O})),\,{\cal
 O})\stackrel{f_q}\longrightarrow {\mathbb C}.
\]
By the isomorphism of {\bf 4.2}:
\begin{equation}
 \otimes_{i=1}^N\Gamma(H^0(D_i,\,\Omega^{\otimes 2}(2{\mathbb O})),\,{\cal
 O})\stackrel{\phi}\simeq {\mathbb C}[\{S_m\}_{m\leq 0}]^{\otimes N},
\end{equation}
it may be considered as a homomorphism from ${\mathbb C}[\{S_m\}_{m\leq
0}]^{\otimes N}$. Now let us recall that each $S_m$ is contained in
the center of the local completion of $U_{-2}(\hat{{\frak G}}_N)$ and
that it acts on $M_{-2,\lambda}$. Hence we can consider a $U_{-2}(\hat{{\frak
G}}_N)$-module $M_{-2,\lambda}^{q}$:
\[
 M_{-2,\lambda}^{q}=M_{-2,\lambda}\otimes_{{\mathbb C}[\{S_m\}_{m\leq
 0}]^{\otimes N}}{\mathbb C}.
\]
Here we have regarded ${\mathbb C}$ as ${\mathbb C}[\{S_m\}_{m\leq
 0}]^{\otimes N}$-module via $f_q$. 
\par\vspace{5pt}
On the other hand, by the homomorphism $\hat{\beta}_{\lambda}$ of {\bf
 Proposition 4.1}, we have a D-module on ${\rm Bun}_{G,X}^{D-fl}$:
\[
 {\cal
 D}^{\prime}_{D-fl,\lambda}\otimes_{\otimes_{i=1}^N\Gamma(H^0(D_i,\,\Omega^{\otimes
 2}(2{\mathbb O})),\,{\cal O})}{\mathbb C}.
\]

\begin{prop}
\[
 \Delta_{D-fl}(M_{-2,\lambda}^{q})\simeq {\cal
 D}^{\prime}_{D-fl,\lambda}\otimes_{\otimes_{i=1}^N\Gamma(H^0(D_i,\,\Omega^{\otimes
 2}(2{\mathbb O})),\,{\cal O})}{\mathbb C}.
\]
\end{prop}
{\bf Proof.}
By definition, we have
\[
 {\cal D}^{\prime}_Z\otimes_{U_{-2}(\hat{{\frak
 G}}_N)}M_{-2,\lambda}^{q}\simeq ({\cal D}^{\prime}_Z\otimes_{U_{-2}(\hat{{\frak
 G}}_N)}M_{-2,\lambda})\otimes_{{\mathbb C}[\{S_m\}_{m\leq
 0}]^{\otimes N}}{\mathbb C}.
\]
The commutativity of $S_m$ with the action of $\hat{B}_+^N$ implies 
\begin{eqnarray*}
 [\tau_{D-fl*}({\cal D}_{Z^{\prime}}\otimes_{U_{-2}(\hat{{\frak G}}_N)}M_{-2,\lambda}^{q})]^{\hat{B}_+^N} &\simeq& [\tau_{D-fl*}({\cal D}_{Z^{\prime}}\otimes_{U_{-2}(\hat{{\frak G}}_N)}M_{-2,\lambda})]^{\hat{B}_+^N}\otimes_{{\mathbb C}[\{S_m\}_{m\leq
 0}]^{\otimes N}}{\mathbb C}\\
&\simeq& {\cal D}^{\prime}_{D-fl,\lambda}\otimes_{{\mathbb C}[\{S_m\}_{m\leq
 0}]^{\otimes N}}{\mathbb C}.
\end{eqnarray*}
Now the isomorphism (17) implies the desired result.
\begin{flushright}
$\triangle$
\end{flushright}
By {\bf Theorem 5.2}, ${\rm Spec}\,D^{\prime}_{D-fl,\lambda}$ may be
identified with $H_{\Delta}(\lambda)$ which is a closed subvariety of
$\oplus_{i=1}^{N}H^0(D_i,\,\Omega^{\otimes 2}(2{\mathbb O}))$ by the
restriction map.
Now the identity 
\[
 \hat{\beta}_{\lambda}(c_0^{(i)})=\Delta(\lambda_i).
\]
and {\bf Proposition 6.1} will imply the following theorem.
\begin{thm}
$\Delta_{D-fl}(M_{-2,\lambda}^{q})$ is nonzero if and only if $q$ is
 global and $\lambda$-admissible.
\end{thm}


\begin{thebibliography}{10}

\bibitem{BNR}
M.~S.~Narasimhan A.~Beauville and S.~Ramanan.
\newblock Spectral curves and the generalized theta divisor.
\newblock {\em J. reine angew. Math.}, 398:169--179, 1989.

\bibitem{Beilinson-Bernstein1993}
A.~Beilinson and J.~Bernstein.
\newblock A proof of {J}antzen conjectures.
\newblock {\em Advances in Soviet Mathematics}, 16(1):1--50, 1993.
\newblock AMS.

\bibitem{Beilinson-Drinfeld}
A.~Beilinson and V.~Drinfeld.
\newblock Quantization of {H}itchin integrable system and {H}ecke eigensheaf.
\newblock Preprint.

\bibitem{Borel}
A.~Borel et~al.
\newblock {\em Algebraic {D}-modules}, volume~2 of {\em Perspective in
  Mathematics}.
\newblock Academic Press, 1987.

\bibitem{Feigin-Frenkel}
B.~Feigin and E.~Frenkel.
\newblock Affine {K}ac-{M}oody algebra at the critical level and
  {G}elfand-{D}ikii algebras.
\newblock {\em Intern. J. of Modern Physics A}, 7:197--215, 1992.

\bibitem{Frenkel1995}
E.~Frenkel.
\newblock Affine algebras, {L}anglands duality and {B}ethe {A}nsatz.
\newblock In {\em Proc. of the Intern. Congr. of Math. and Phys.}, pages
  606--642. International Press, 1995.

\bibitem{Frenkel02}
E.~Frenkel.
\newblock Lectures on {W}akimoto modules, opers and the center at the critical
  level.
\newblock math.QA/0210029, November 2002.

\bibitem{Frenkel2003}
E.~Frenkel.
\newblock Recent advances in the {L}anglands program.
\newblock math.AG/0303074, March 2003.

\bibitem{Frenkel-BenZvi}
E.~Frenkel and D.~Ben-Zvi.
\newblock {\em Vertex {A}lgebras and {A}lgebraic {C}urves}, volume~88 of {\em
  Mathematical Surveys and Monographs}.
\newblock AMS, 2001.

\bibitem{Ginzburg}
V.~Ginzburg.
\newblock Perverse sheaves on a loop group and {L}anglands' duality.
\newblock alg-geom/9511007, November 2000.

\bibitem{Hitchin1987}
N.~Hitchin.
\newblock Stable bundles and integrable systems.
\newblock {\em Duke Math. J.}, 54(1):91--114, 1987.

\bibitem{LS}
Y.~Laszlo and C.~Sorger.
\newblock The line bundles on the moduli of parabolic {G}-bundles over curves
  and their sections.
\newblock {\em Ann. Sci. Ecole Norm. Sup.}, 30:499--525, 1997.

\bibitem{Laumon}
G.~Laumon.
\newblock Un analogue global du c\^{o}ne nilpotent.
\newblock {\em Duke Math. Jour.}, 57:647--671, 1988.

\bibitem{Laumon-MoretBailly}
G.~Laumon and L.~Moret-Bailly.
\newblock {\em Champs alg\'{e}briques}, volume~39 of {\em A series of Modern
  Surveys in Mathematics}.
\newblock Springer, 1999.

\bibitem{Sorger}
C.~Sorger.
\newblock Lectures on moduli of principal {G}-bundles over algebraic curves.
\newblock Lectures given at School on Algebraic Geometry Trieste, 26 July -13
  August 1999.

\end{thebibliography}

\vspace{10mm}
\begin{flushright}
Address : Department of Mathematics and Informatics\\
Faculty of Science\\
Chiba University\\
1-33 Yayoi-cho Inage-ku\\
Chiba 263-8522, Japan \\
e-mail address : sugiyama@math.s.chiba-u.ac.jp
\end{flushright}

\end{document}